\documentclass[12pt]{article}
\usepackage{amsfonts}
\usepackage{amssymb}
\newtheorem{lemma}{Lemma}

\newtheorem{proposition}[lemma]{Proposition}
\newtheorem{corollary}[lemma]{Corollary}

\def\pf{\mbox{\bf Proof. }}

\begin{document}
\title{Some remarks on varieties of pairs of
commuting upper triangular matrices and an interpretation of
commuting varieties}
\author{\small{Roberta Basili}\\ \small{Via dei Ciclamini 2 B 06126 Perugia Italy, e-mail:
robasili@alice.it}}
\date{}
\maketitle \begin{abstract} It is known that the variety of pairs
of $n\times n$ commuting upper triangular matrices isn't a
complete intersection for infinitely many values of $n$; we show
that there exists $m$ such that this happens if and only if $n>m$.
We also show that $m<18$ and that it could be found by determining
the dimension of the variety of pairs of commuting strictly upper
triangular matrices. Then we define a natural map from the variety
of pairs of commuting $n\times n$ matrices onto a subvariety
defined by linear equations of the grassmannian of subspaces of
$K^{n^2}$ of codimension 2.\end{abstract}
\section{Introduction}Let $T_n$ be the set of all $n\times n$ upper triangular matrices
over an algebraically closed field $K$; let ${\cal T}_n$ be the
subset of $T_n$ of all the invertible matrices. Let
$$ CT_n=\{ (X,Y)\in T_n\times T_n\, :\, [X,Y]=0\}\, .$$
Let $U_n$ be the subset of $T_n$ of all the strictly upper
triangular matrices and let
$$NT_n=CT_n\cap \big( U_n\times U_n\big) \, .$$
It is known that there exist infinitely many values of $n$ such
that $CT_n$ and $NT_n$ are not irreducible and are not complete
intersections. The determination of the smallest $n$ such that
these properties occur is an open problem which has recently
interested several mathematicians. \newline The action of ${\cal
T}_n$ on $U_n$ hasn't finitely many orbits; a classification of
them can be found in \cite{He}. Hence many of the arguments which
are used in the study of commuting varieties cannot be applied
here.\newline In section \ref{s2} we show that $CT_n$ is a
complete intersection if and only if its irreducible components
have the same dimension, and that there exist natural numbers
$m,\; m'$ such that $CT_n$ isn't a complete intersection if and
only if $n>m$ and is reducible if and only if $n>m'$. Similar
results hold for $NT_n$, moreover we prove that $m$ and $m'$ have
the previous property according to the dimension of $NT_m$,
$NT_{m'}$. Then we give examples which prove that $m<18$ and
$m'<17$.
\newline Many results of section \ref{s2} were independently
obtained by Allan Keeton, as results of his Ph.D. thesis. More
precisely, Keeton communicated to the author that he had obtained
the following results:
\begin{itemize} \item[{\em a)}] $CT_n$ is an irreducible, normal
complete intersection if $n\leq 8$; \item[{\em b)}] $CT_n$ is not
normal if $n\geq 16$; \item[{\em c)}] $CT_n$ is reducible if
$n\geq 17$; \item[{\em d)}] $CT_n$ is not a complete intersection
and not of pure dimension if $n\geq 18$.\end{itemize} In section
\ref{s3} we consider the variety ${\cal  C}(n,K)$ of all the pairs
of commuting $n\times n$ matrices, which we regard as a subvariety
of $\mathbb P^{n^2-1}\times \mathbb P^{n^2-1}$. We denote by
${\cal C}_0(n,K)$ the subvariety of ${\cal C}(n,K)$ of pairs of
equal elements, then we define a map $\gamma _n$ from ${\cal
C}(n,K)\setminus {\cal C}_0(n,K)$ into the grassmannian $G(2,
K^{n^2})$ of all the subspaces of $K^{n^2}$ of codimension 2. The
fibers of this map are the orbits of ${\cal C}(n,K)\setminus {\cal
C}_0(n,K)$ under the natural action of GL$(2,K)$ on ${\cal
C}(n,K)\setminus {\cal C}_0(n,K)$. We get that the image of
$\gamma _2$ is a linear complete intersection subvariety of the
projective space of dimension 5 in which $G(2,K^4)$ is defined.
\section{Some remarks on varieties of pairs of commuting upper triangular
matrices}\label{s2} We will denote by $(X,Y)$ a generic element of
$T_n\times T_n$. The entries of $[X,Y]$ give $\ {\displaystyle
\frac{n(n-1)}{2}}\ $ equations for $CT_n$ and ${\displaystyle
\frac{(n-1)(n-2)}{2}}$ equations for $NT_n$.\newline Let $CT_n^0$
be the Zariski closure of the subset of $CT_n$ of all the pairs
$(X,Y)$ such that $X$ and $Y$ are regular (that is have minimum
polynomial of degree $n$); let $NT_n^0$ be the same Zariski
closure in $NT_n$.\begin{proposition} We have: \begin{itemize}
\item[i)] $CT_n^0$ is irreducible of dimension ${\displaystyle
\frac{n(n+3)}{2}}$; this is the minimum dimension of the
irreducible components of $CT_n$. \item[ii)] $NT_n^0$ is
irreducible of dimension ${\displaystyle \frac{n(n+1)}{2}-1}$;
this is the minimum dimension of the irreducible components of
$NT_n$.\end{itemize}\end{proposition} \pf If we consider the
canonical projection of $CT_n^0$ on $T_n$, the fibers of the
regular matrices have dimension $n$, hence $CT_n^0$ is irreducible
of dimension $${\displaystyle \frac{n(n+1)}{2}+n}={\displaystyle
\frac{n(n+3)}{2}}\, .$$ Moreover
$${\displaystyle \frac{n(n+3)}{2}}=\mbox{\rm dim }\big( T_n\times T_n\big) -
{\displaystyle \frac{n(n-1)}{2}}\, ,$$ which shows i). The same
argument can be used for ii).\vspace{3mm}\newline By the
irreducibility of the centralizer in $T_n$ and in $U_n$ we get the
following result.
\begin{proposition} If $\big( X,Y\big) \in CT_n $ \big( $NT_n $
\big) and $\, X$ or $Y$ commutes with regular matrices of $T_n$
\big( $U_n$ \big) then $(X,Y)\in CT_n^0$ \big( $NT_n^0$\big)
.\end{proposition} We observe that any irreducible component of
$CT_n$ \big( $NT_n$\big) is stable under the action of ${\cal
T}_n$. Moreover any irreducible component of $CT_n$ is stable with
respect to the action of $K^2$ defined by
$$(x,y)\cdot (X,Y)=(X+xI_n,Y+yI_n)\, ;$$
hence the subset of the pairs of nonsingular matrices is dense in
any irreducible component of $CT_n$.\vspace{3mm}\newline We denote
by $\{ e_1,\ldots ,e_n\} $ the canonical basis of $K^n$ and by
$M(p,q)$ the set of all $p\times q$ matrices over $K$.
\begin{proposition} \label{3} If $CT_{n-1}$ \big( $NT_{n-1}$\big) isn't
 irreducible or
isn't a complete intersection, the same holds for $CT_{n}$ \big(
$NT_{n}$\big).\end{proposition}  \pf We first prove the claim for
$CT_{n}$.\newline Let $CT_{n-1}^1$ be an irreducible subvariety of
$CT_{n-1}$ different from $CT_{n-1}^0$ and let $$\mbox{\rm dim
}CT_{n-1}^1={\displaystyle \frac{(n-1)(n+2)}{2}}+k\, ,\qquad k\geq
0\, .$$ Let $T_{n}'$ be the subspace of $T_{n}$ of all the
endomorphisms which stabilize $\langle e_1\rangle $ and $\langle
e_2,\ldots ,e_{n}\rangle $. Let ${\cal T}_{n}'={\cal T}_{n}\cap
T_{n}'$. Let
$$\Gamma =\left\{ (X,Y)\in CT_{n}\, : \right. $$ $$\left.
X=\pmatrix{ 0 & 0\cr 0 & X'\cr} \, , \ Y=\pmatrix{ 0 & 0\cr 0 &
Y'\cr}\, , (X',Y')\in CT_{n-1}^1\, \right\} \, .$$ Let $\Gamma '$
be the orbit of $\Gamma $ under the action of ${\cal T}_{n}$. If
$(X',Y')\in CT_{n-1}^1$, rank $X'$, rank $Y'$ $=n-1$ and $G\in
{\cal T}_{n}$ then we have that $G\cdot (X,Y)\in \Gamma \, $ iff
$\, G\in {\cal T}_{n}'$. Hence
$$\mbox{\rm dim } \Gamma '=\mbox{\rm dim }{\cal T}_{n}
-\mbox{\rm dim }{\cal T}_{n}'+\mbox{\rm dim }CT_{n-1}^1$$ $$=
{\displaystyle \frac{n(n+1)}{2}}-\bigg({\displaystyle
\frac{n(n-1)}{2}}+1\bigg)+{\displaystyle
\frac{(n-1)(n+2)}{2}}+k=$$
$$={\displaystyle \frac{n(n+3)}{2}}+k-2\, .$$
Let $CT_{n}'$ be the orbit of $\Gamma '$ under the action of
$K^2$; we have $$\mbox{\rm dim }CT_{n}^1={\displaystyle
\frac{n(n+3)}{2}}+k\, .$$ Since $X'$ and $Y'$ aren't regular,
$CT_{n}^1\neq CT_{n}^0$, which shows the claim. \newline We now
prove the claim for $NT_{n}$. Let $NT_{n-1}^1$ be a subvariety of
$NT_{n-1}$ and let
$$NT_{n}^1 =\left\{ (X,Y)\in NT_{n}\, :\,
X=\pmatrix{ 0 & \widetilde X\cr 0 & X'\cr}\; ,\right. $$ $$\left.
Y=\pmatrix{ 0 & \widetilde Y\cr 0 & Y'\cr}\, , (X',Y')\in
NT_{n-1}^1\, , \widetilde X,\widetilde Y\in M(1,n-1)\right\} \,
.$$  The equations for $NT_{n}^1$ as subvariety of
$$NT_{n-1}^1\times M(1,n-1)\times M(1,n-1)$$ are given by
$\widetilde XY'-\widetilde YX'=0$, hence
$$\mbox{\rm dim }NT_{n}^1\geq \mbox{\rm dim }NT_{n-1}^1+n.$$ If
$\mbox{\rm dim }NT_{n-1}^1={\displaystyle \frac{n(n-1)}{2}}-1+k\,
,\ k\geq 0$ then $\mbox{\rm dim
 }NT_{n}\geq {\displaystyle
\frac{n(n+1)}{2}}-1+k\, ,$ which proves the
claim.\vspace{3mm}\newline For $X\in T_n$ let $f(t)=(t-\lambda
_1)^{m_1}\cdot \cdots \cdot (t-\lambda _r)^{m_r}$ be the minimum
polynomial of $X$. For $i=1,\ldots ,r$ let $f_i(t)=f(t):(t-\lambda
_i)^{m_i}$ and let $g_1(t),\ldots ,g_r(t)$ be such that $\,
{\displaystyle \sum _{i=1}^rg_i(t)f_i(t)=1}$; then $\,
{\displaystyle \sum_{i=1}^rg_i(X)f_i(X)=I_n}\, $. Hence the
matrices $I_i=g_i(X)f_i(X)\, $, $i=1,\ldots ,r$, are orthogonal
projections of $K^n$ on $K^n$ and the image of $I_i$ is ker
$(X-\lambda _iI_n)^{m_i}$. Then we get the following result.
\begin{lemma} \label{4} Let $(X,Y)\in CT_n$. There exist $G\in {\cal T}_n$
and a partition $\{ E_1,\ldots ,E_r\} \ $ of $\ \{ Ge_1,\ldots
,Ge_n\} \ $ such that $\mbox{\rm ker }(X-\lambda
_iI_n)^{m_i}=\langle E_i\rangle $ and $\langle E_i\rangle $ is
stable with respect to $X$ and $Y$ for $i=1,\ldots ,r$.\end{lemma}
\pf The matrices $I_i$ for $i=1,\ldots ,r$ are upper triangular
and commute with the matrices of the centralizer of $X$. If $j\in
\{1,\ldots ,n\} $ there exists a unique $i\in \{1,\ldots ,r\} $
such that the entry of $I_i$ of indices $(j,j)$ is 1. Let $G\in
{\cal T}_n$ be such that $Ge_j=I_ie_j$; this gives a partition
with the required property.\vspace{3mm}\newline We set $n_i=|E_i|$
for $i=1,\ldots ,r$. By Lemma \ref{4} we get the following
results.
\begin{proposition} \label{5} Let $CT_{n-1}$ be a complete intersection and
let $(X,Y)\in CT_n$ be such that $X$ or $Y$ has at least two
eigenvalues. Then $(X,Y)$ doesn't belong to any irreducible
component of dimension greater than ${\displaystyle
\frac{n(n+3)}{2}}$.\end{proposition} \pf If $(X,Y)$ belongs to an
irreducible component then the subset of it of all the pairs such
that at least one of the matrices has more than one eigenvalue is
dense. Let $E=\{ E_1,\ldots ,E_r\} $ be a partition of
$\{e_1,\ldots ,e_n\}$ such that $r\geq 2$. Let $T_E$ be the subset
of $T_n$ of all the endomorphisms which stabilize $\langle
E_i\rangle $ for $i=1,\ldots ,r$. Let ${\cal T}_E ={\cal T }_n\cap
T_E$ and let $CT_E=CT_n\cap (T_E\times T_E)$. By Proposition
\ref{3} we have dim $CT_E={\displaystyle \sum_{i=1}^r
{\displaystyle \frac{n_i(n_i+3)}{2}}}$, hence the dimension of the
orbit of $CT_E$ under the action of ${\cal T}_n$ is less or equal
than
$$\mbox{\rm dim }{\cal T}_n-\mbox{\rm dim }{\cal T}_E+\mbox{\rm
dim }CT_E= $$ $$ ={\displaystyle
\frac{n(n+1)}{2}}-\sum_{i=1}^r{\displaystyle
\frac{n_i(n_i+1)}{2}}+\sum_{i=1}^r{\displaystyle
\frac{n_i(n_i+3)}{2}}={\displaystyle \frac{n(n+3)}{2}}\, .$$ Hence
the claim follows by Lemma \ref{4}.
\begin{proposition} \label{6} Let $CT_{n-1}$ be irreducible and
let $(X,Y)\in CT_n$ be such that $X$ or $Y$ has at least two
eigenvalues. Then $(X,Y)\in CT_n^0$. \end{proposition} \pf Let us
assume that $X$ has $r\geq 2$ eigenvalues and let $\{ E_1,\ldots
,E_r\} $ be as in Lemma \ref{4}. Let $CT_n^E$ be the subvariety of
$CT_n$ of all the pairs of matrices which stabilize $\langle
E_i\rangle $ for $i=1,\ldots ,r$. Then we have
$$CT_n^E\cong CT_{n_1}\times \cdots \times CT_{n_r}\, .$$
By Proposition \ref{3} $\ CT_{n_i}$ is irreducible for $i=1,\ldots
,r$. Hence $CT_n^E$ is irreducible and the subset of $CT_E$ of all
the pairs of regular matrices is dense, which shows the claim.
\vspace{3mm}\newline By Propositions \ref{5} and \ref{6} we get
that, in order to determine the values of $n$ such that $CT_n$
isn't irreducible or isn't a complete intersection, we can look
for irreducible components which have as elements only pairs of
matrices with only one eigenvalue. Hence it is enough to determine
the dimension of $NT_n$.\newline Let $X=(x_{i,j})$, $Y=(y_{i,j})$,
$X_{i,j}=\pmatrix{x_{i,j} \cr y_{i,j}\cr }$. The condition
$[X,Y]=0$ gives the following equations for $NT_n$:
\begin{equation}\label{e} \sum_{k=i+1}^{j-1} \mbox{\rm det }\pmatrix{ X_{i,k} &
X_{k,j}}=0\quad  i= 1,\ldots ,n-2\, ,\ j=i+2,\ldots, n\,
.\end{equation} We observe that this system of equations is
invariant under the involution of $U_n$ defined by $z_{i,j}\mapsto
z_{n+1-j,n+1-i}$.\newline Any irreducible component of $NT_n$
different from $NT_n^0$ is contained in a subvariety of $NT_n$
defined by some equations of the form $X_{i,j}=0$, $i<j$. Moreover
the following result holds.
\begin{lemma} \label{7} If $X\in U_n$ and there exist
$h,k\in \{ 1,\ldots ,n\} $ such that $x_{h,k}$ vanishes on the
orbit of $X$ under the action of ${\cal T}_n$ then also $x_{i,j}$
vanishes on that orbit for $i,j=h,\ldots ,k$.\end{lemma} \pf The
claim follows from the fact that in the orbit under ${\cal T}_m$
of any nonzero matrix of $U_m$ there exists a matrix such that its
entry of indices $(1,m)$ isn't 0.
\begin{corollary} \label{8} Let $NT_n^{\ast }$ be an irreducible
component of $NT_n$ different from $NT_n^0$. There exists $s\in
\{1,\ldots ,n-1\} $ and subsets $\ J_1,\ldots ,J_s\ $ of $\ \{
1,\ldots, n\} \ $,  such that: \begin{itemize} \item[i)]  if $h\in
\{ 1,\ldots ,s\} \, , $ $i,j\in J_h$, $l\in \{1,\ldots,n\} $ and
$i<l<j$ then $l\in J_h$; \item[ii)] $J_1\cup \cdots \cup J_s =\{
1,\ldots ,n\} $;  \item[iii)] $X_{i,j}$ is $0$ on $NT_n^{\ast }$
iff there exists $h\in \{1,\ldots ,s\} $ such that $i,j\in J_h$.
 \end{itemize} \end{corollary}
 Let $\Upsilon _n$ be the set of all the partitions $J= \{ J_1,\ldots ,J_s\} \
 $ of $\ \{ 1,\ldots
 ,n\} $ such that $s\in \{ 1,\ldots ,n-1\} $ and $J_1,\ldots ,J_s$ have the property i) of
 Corollary \ref{8}. We
assume that if $h,k\in \{1,\ldots ,s\} $ and $h<k$ the elements of
$J_h$ are smaller than those of $J_k$. If $J\in \Upsilon _n$ we
denote by $U_n^J$ the subvariety of $U_n$
 defined by the equations $$x_{i,j}=0\, ,\qquad i,j\in J_h\, ,\ h=
 1,\ldots ,s \; $$ and we set $NT_n^J=NT_n\cap (U_n^J\times
 U_n^J)$.
\newline If $NT_n^{\ast }$ is an irreducible
 component of $NT_n$ different from $NT_n^0$ there exists $J\in
 \Upsilon _n$ such that $NT_n^{\ast }\subseteq
 NT_n^J$.\vspace{3mm}\newline
{\bf Example A}\hspace{2mm} The variety $NT_4$ is defined by the
equations:
$$\mbox{\rm det }\pmatrix{X_{1,2} & X_{2,3}\cr}=0\, ,\qquad
\mbox{\rm det }\pmatrix{X_{2,3} & X_{3,4}\cr}=0\, ,$$ $$ \mbox{\rm
det }\pmatrix{X_{1,2} & X_{2,4}\cr}+\mbox{\rm det
}\pmatrix{X_{1,3} & X_{3,4}\cr}=0\; .$$ Let $NT_4^1$ be the
subvariety of $NT_4$ defined by the equations $X_{2,3}=0$ and let
$NT_4^0$ be the subvariey of $NT_n$ defined by the equations
$$ \mbox{\rm rank }\pmatrix{X_{1,2} & X_{2,3} & X_{3,4}\cr}\leq 1\,
.$$ We have that $NT_4^1$ and $NT_4^0$ are proper subvarieties of
$NT_4$ such that $NT_4=NT_4^1\cup NT_4^0$. It is obvious that
$NT_4^1$ is irreducible; let us prove that $NT_4^0$ is the closure
of the subset of $NT_4$ of all the pairs of regular matrices.
\newline Let ${\cal A}$ be an open subset of $NT_4^0$ and let
$(\overline X,\overline Y)\in {\cal A}$. We may assume that
$\overline X_{2,3}\neq 0$. Let $\overline X_{1,2}, \overline
X_{3,4}=0$. The subspace of $U_4\times U_4$ given by the equations
$X_{1,2}=0$ and $X_{3,4}=0$ is contained in $NT_4^0$, hence we can
assume that $\overline X_{2,3}$, $\overline X_{2,4}$, $\overline
X_{1,3}$ are pairwise linearly independent. Let $\alpha, \beta \in
K$ be such that $\overline X_{2,3}=\alpha \overline X_{2,4}+\beta
\overline X_{1,3}$; let us consider the subvariety of $NT_n^0$,
parametrized by $\tau $, defined by: $X_{i,j}=\overline X_{i,j}$
for $(i,j)=(2,3), (2,4), (1,3)$, $X_{1,2}=\tau \alpha \overline
X_{2,3}$, $X_{3,4}=- \tau \beta \overline X_{2,3}$; for $\tau \neq
0$ we get pairs of regular elements, which shows the claim. Let
$\overline X_{1,2}=\gamma \overline X_{2,3}$, $\gamma \neq 0$, and
let $\overline X_{3,4}=0$. Let us consider the subvariety of
$NT_n^0$, parametrized by $\delta $, defined by $X_{i,j}=\overline
X_{i,j}$ for $(i,j)=(2,3),(1,2),(1,3)$, $X_{3,4}=\delta \overline
X_{2,3}$, $X_{2,4}=\overline X_{2,4}+{\displaystyle \frac{\delta
}{\gamma }}\overline X_{1,3}$; for $\delta \neq 0$ we get pairs of
regular elements. We could use a similar argument if $\overline
X_{1,2}=0$ and $\overline X_{3,4}\neq 0$, which shows the
claim.\vspace{3mm}\newline Let $J\in \Upsilon _n$ and let
$Z_{h,k}=\big( x_{i,j}\big)\ $, $\ W_{h,k}=\big( y_{i,j}\big)$, $\
i\in J_h,\ j\in J_k$. We can write the equations of $NT_n^J$ in
$U_n^J\times U_n^J$ as follows:
$$\sum_{i=h+1}^{k-1}Z_{h,i}W_{i,k}-W_{h,i}Z_{i,k}=0\qquad h=1,\ldots ,s-2,\ \;
k=h+2,\ldots ,s \;.$$ This can be also written in the following
way:
$$\pmatrix{ Z_{h,h+1} & -W_{h,h+1} & \cdots & Z_{h,k-1} & -W_{h,k-1}\cr
}\; \pmatrix{ W_{h+1,k}\cr Z_{h+1,k}\cr \vdots \cr W_{k-1,k}\cr
Z_{k-1,k}\cr}\ =\ 0 $$ $$h=1,\ldots ,s-2,\ \; k=h+2,\ldots ,s
\;.$$ Let $V_{m,p,q}=\{ (A,B)\in M(m,p)\times M(p,q)\, :\, AB=0\}
$. We can determine a lower bound of the dimension of $NT_n^J$ by
the following elementary result.
\begin{lemma}\label{11} The irreducible components of $V_{m,p,q}$
are the subvarieties
$$V_{m,p,q}^{a,b}=\{ (A,B)\in V_{m,p,q}\, :\, \mbox{\rm rank
}A\leq a,\ \mbox{\rm rank }B\leq b\} $$ where $(a,b)$ is maximal
such that $\; b\leq \mbox{\rm min } \{ p,q \} \; $,  $\; a\leq
\mbox{\rm min } \{ p-b,m \}\; .$ We have:
\begin{itemize} \item[i)] $\mbox{\rm dim
}V_{m,p,q}^{a,b}=a(p+m-a)+b(p+q-b)-ab\; ;$ \item[ii)] $V_{m,p,q}$
is a complete intersection iff $\; p\geq m+q-1$.\end{itemize}
\end{lemma} By the following example W.V. Vasconcelos observed
that $CT_n$ isn't a complete intersection for infinitely many
values of $n$.\vspace{3mm}\newline {\bf Example B} \cite{Va}. Let
$n=3m$ and let $J\in \Upsilon _{3m}$ be defined by $J_1=\{
1,\ldots ,m\} $, $J_2=\{ m+1,\ldots, 2m\} $, $J_3=\{ 2m+1,\ldots,
3m\} $. If $(X,Y)$ belongs to the subvariety of $U_n\times U_n$
defined by
 $$X_{i,j}=0\, ,\qquad i,j\in J_h,\quad h=1,2,3$$
 then for $i\not \in J_3$ or $j\not \in J_3$ the entry of $[X,Y]$
 of indices $(i,j)$
is $0$. Hence $\mbox{\rm dim }NT_{3m}^J\geq 3m^2+3m^2-m^2=5m^2$.
We have dim $CT_{3m}^0={\displaystyle \frac{9m(m+1)}{2}}$, which
for $m\geq 10$ is smaller than $5m^2$, hence $CT_{3m}$ isn't a
complete intersection for $m\geq 10$.\vspace{3mm}\newline {\bf
Example C}\hspace{2mm} Let $n=18$ and let $J\in \Upsilon _{18}$ be
such that $|J_1|=1$, $|J_2|=5$, $|J_3|=6$, $|J_4|=5$, $|J_5|=1$.
The condition $[X,Y]=0$ gives 48 equations for $NT_{18}^J$ as
subvariety of $U_{18}^J\times U_{18}^J$. Hence
$$\mbox{\rm dim } NT_{18}^J\geq \mbox{\rm dim }(U_{18}^J\times
U_{18}^J)-48=188\; .$$ Then the dimension of the orbit of
$NT_{18}^J$ under the action of $K^2$ is greater or equal than
$190$. Since $\mbox{\rm dim }CT_{18}^0=189$, $CT_{18}$ isn't a
complete intersection.\vspace{3mm}\newline {\bf Example
D}\hspace{2mm} Let $n=17$ and let $J\in \Upsilon _{1}$ be such
that $|J_1|=2$, $|J_2|=4$, $|J_3|=5$, $|J_4|=4$, $|J_5|=2$. The
condition $[X,Y]=0$ gives 56 equations for $NT_{17}^J$ as
subvariety of $U_{17}^J\times U_{17}^J$. Hence
$$\mbox{\rm dim } NT_{17}^J\geq \mbox{\rm dim }(U_{17}^J\times
U_{17}^J)-56=168\; .$$ Then the dimension of the orbit of
$NT_{17}^J$ under the action of $K^2$ is greater or equal than
$170$. Since $\mbox{\rm dim }CT_{17}^0=170$, $CT_{17}$ is
reducible. \section{An interpretation of commuting
varieties}\label{s3} Let $M(n,K)$ be the set of $n\times n$
matrices over K, which we regard as a projective space of
dimension $n^2-1$. Let
$${\cal C}(n,K)=\{ (X,Y)\in M(n,K)\times M(n,K)\, :\, [X,Y]=0\}
\; .$$ For $X,Y\in M(n,K)$ we set $X=\big( x_{i,j}\big) $,
$Y=\big( y_{i,j}\big) $, where $i,j\in \{1,\ldots ,n\} $, and
$X_{i,j}=\pmatrix{ x_{i,j} \cr y_{i,j}\cr }$. As a generalization
of equations \ref{e}, the equations for ${\cal C}(n,K)$ given by
the condition $[X,Y]=0$ can be written as follows:
$$ \sum_{k=1}^{n} \mbox{\rm det }\pmatrix{ X_{i,k} &
X_{k,j}}=0\qquad  i,\, j=1,\ldots, n\, .$$ Let $${\cal
C}_0(n,K)=\{ (X,Y)\in {\cal C}(n,K)\, :\, \mbox{\rm det }\pmatrix{
X_{i,j} & X_{h,k}\cr }=0 \ \ \  \mbox{\rm for any }$$ $$
\hspace{15mm}(i,j),\ (h,k)\in \{ 1,\ldots ,n\} \times \{ 1,\ldots
, n\} \} \, .$$ For $(i,j),\ (h,k)\in \{ 1,\ldots ,n\} \times \{
1,\ldots ,n\} $ we denote by $p_{(i,j)(h,k)}$ the Pl\"{u}cker
coordinates of subspaces of codimension $2$ of $K^{n^2}$; we
denote by $G(2,K^{n^2})$ the grassmannian of those subspaces.
\newline There is a natural map $\gamma _n$ from ${\cal
C}(n,K)\setminus {\cal C}_0(n,K)$ into $G(2,K^{n^2})$, defined by
associating to $(X,Y)$ the subspace having the following
Pl\"{u}cker coordinates:
$$p_{(i,j)(h,k)}=\mbox{\rm det }\pmatrix{ X_{i,j} & X_{h,k}}$$ for
$(i,j),\ (h,k)\in \{ 1,\ldots ,n\} \times \{ 1,\ldots ,n\}$. The
image of $\gamma _n$ is the subvariety of $G(2,K^{n^2})$ defined
by the following linear equations:
$$\sum _{k=1}^n p_{(i,k)(k,j)}=0 \ , \hspace{30mm} i,j=1,\ldots
,n\; .$$  The group GL$(2,K)$ acts on ${\cal C}(n,K)\setminus
{\cal C}_0(n,K)$ by the following rule:
$$\pmatrix{a & b \cr c & d\cr }\; \cdot \;
(X,Y)=(a X+b Y, cX+dY)$$ for $\pmatrix{a & b \cr c & d\cr }\in
\mbox{\rm GL}(2,K)$ and $(X,Y)\in {\cal C}(n,K)\setminus {\cal
C}_0(n,K)$. The fibers of $\gamma _n$ are the orbits of ${\cal
C}(n,K)\setminus {\cal C}_0(n,K)$ under the action of
GL$(2,K)$.\newline We can define a similar map for $CT_n$ and
$NT_n$, as a restriction of $\gamma _n$. As an example we
illustrate this geometrical interpretation for
$NT_4$.\vspace{3mm}\newline {\bf Example A} \hspace{2mm} We regard
$U_4$ as a projective variety of dimension 5 (whose coordinates
have indices $(1,2)$, $(2,3)$, $(3,4)$, $(1,3)$, $(2,4)$,
$(1,4)$). We consider the elements of $U_4$ as hyperplanes of
$\mathbb P^5$. The map $\gamma _4$ associates to any pair of
different hyperplanes of $NT_4$ the subspace of $\mathbb P^5$
given by their intersection. The image of $NT_4$ by $\gamma _4$ is
defined by the equations:
$$p_{(1,2)(2,3)}=0\; ,\ p_{(2,3)(3,4)}=0\; ,\ p_{(1,2)(2,4)}+
p_{(1,3)(3,4)}=0\; .$$ The inverse image under $\gamma _4$ of the
subset of all subspaces of $\mathbb P^5$ of codimension 2 such
that $p_{(1,2)(3,4)}=0$ is the irreducible component $NT_4^0$ of
$NT_4$. The set of pairs of hyperplanes such that the coordinate
of indices $(2,3)$ is $0$ (that is, "parallel" to the
$(2,3)$-axis) is the irreducible component $NT_4^1$ of $NT_4$.
\vspace{3mm}\newline {\bf Example E}\hspace{2mm} The image of
${\cal C}(2,K)\setminus {\cal C}_0(2,K)$ under the map $\gamma _2$
is a subvariety of $G(2,K^{4})$ defined by the following
equations:
\begin{equation}\label{e2} p_{(1,2)(2,1)}=0\; ,\end{equation}
\begin{equation}\label{e3} p_{(1,1)(1,2)}+p_{(1,2)(2,2)}=0\;
,\end{equation} \begin{equation}\label{e4}
p_{(2,1)(1,1)}+p_{(2,2)(2,1)}=0\; .\end{equation} The variety
$G(2,K^{4})$ is a subvariety of a projective space of dimension 5
defined by the equation \begin{equation}\label{e5}
p_{(1,1)(2,2)}p_{(1,2)(2,1)}-p_{(1,1)(1,2)}p_{(2,2)(2,1)}+p_{(1,1)(2,1)}p_{(2,2)(1,2)}=0\;
.\end{equation} If we consider subvarieties of that projective
space, we have the following results. The subvariety defined by
equations \ref{e3}, \ref{e4} and \ref{e5} has two irreducible
components; one of them is $\gamma_2({\cal C}(2,K)\setminus {\cal
C}_0(2,K))$ (see \cite{Kn}). The subvariety defined by the
equations \ref{e2}, \ref{e3} and \ref{e4} is contained in $G(2,
K^4)$, hence it is $\gamma_2({\cal C}(2,K)\setminus {\cal
C}_0(2,K))$. Then $\gamma_2({\cal C}(2,K)\setminus {\cal
C}_0(2,K))$ is linear complete intersection as subvariety of that
projective space. \vspace{6mm}\newline{\bf
Acknowledgement}\hspace{2mm} The problem of properties of
commuting varieties of upper triangular matrices was suggested to
the author by Prof. Jerzy Weyman during her stay at Northeastern
University in summer 2003.

\end{document}